\documentclass[11pt]{article}
\usepackage{graphicx}
\usepackage{amsmath}
\usepackage{amsfonts}
\usepackage[english]{babel}
\usepackage{amsthm}
\newtheorem{theorem}{Theorem}[section]
\newtheorem{proposition}{Proposition}[section]

\theoremstyle{definition}

\theoremstyle{remark}
\newtheorem{remark}{Remark}[section]
\theoremstyle{example}

\begin{document}

\title{On a Side Condition for \\ Wronskian-Involving Differential Equations}
\author{Nicoleta B\^{\i}l\u{a}
\footnote{Department of Mathematics and Computer Science, Fayetteville State University,
1200 Murchison Road, Fayetteville, NC 28301, E-mail: \texttt{nbila@uncfsu.edu}.}}
\date{}
\maketitle

\begin{abstract}
The purpose of this paper is to make a few connections among specific concepts occurring in differential geometry and the theory of differential equations with the aim of identifying an intriguing class of undetermined nonlinear ordinary differential equations whose solutions satisfy a specific side condition consisting in a homogeneous third-order linear ordinary differential equation. A method for solving this class of Wronskian-involving differential equations based on the proposed side condition is presented. The Tzitzeica curve equation arising in the theory of space curves is considered as an example, and new closed and integral-form solutions for this equation are obtained.
\end{abstract}


\section{Introduction}\label{Section1}

Since there is no general theory for integrating nonlinear ordinary differential equations (ODEs),
it will always be a challenge to explore methods involving intriguing side conditions that would allow one ``to undo" the nonlinearity encapsulated among the unknown functions and their derivatives and obtain particular solutions to the studied equation. The classical Lie method \cite{OlverBook} is one of the most well-known techniques that may be applied in this situation and, for instance, would allow one to reduce the order of an ODE by one whenever the equation is invariant with respect to a particular one-parameter Lie group of transformations.
However, this reduction may become more complicated in the case of underdetermined nonlinear ODEs of higher order that have fewer equations than unknowns  and contain arbitrary constants or functions. In this paper, the focus is on a specific class of nonlinear differential equations (\ref{EqnWronskian}) that expresses a relation between the Wronskian (\ref{Wronskian1}) of the unknown functions and the Wronskian (\ref{WronskianDer}) of their first derivatives. This class of autonomous ODEs is analyzed in the particular case when its solutions satisfy the auxiliary equation (\ref{EqnU3}). For instance, the Tzitzeica curve equation (\ref{TzitzeicaWronskian}) belongs to this class. It is important to point out that the side condition (\ref{EqnU3}) is not found by using the classical Lie method but rather is related to how the nonlinearity among the unknown functions and their derivatives has occurred while deriving the equation.

Introduced in 1812 by J\'{o}zef Maria Ho\"{e}n\'{e}-Wro\'{n}ski and later on mentioned by Thomas Muir in \cite{Muir}, the \textit{Wronskian determinant} (or, simply, the \textit{Wronskian}) of a set of $n$ smooth functions is the determinant of the $n\times n$ matrix whose entries are the functions and their derivatives up to the $(n-1)$st order. Namely, the functions are listed in the first row, their first derivatives are placed on the second row, and so on, and their $(n-1)$st derivatives are written on the last row (\cite{BoyceDiPrima}, p. 221). In particular, for $n=3$, the Wronskian $W(x,y,z)(t)$ of the smooth functions $x(t)$, $y(t)$, and $z(t)$ is given by (\ref{Wronskian1}). If the functions represent the set of fundamental solutions of a $n$th-order homogeneous linear ODE, then their Wronskian does not vanish and satisfies the Abel-Liouville-Ostrogradski formula (see, e.g., \cite{BoyceDiPrima}, p. 239). Interestingly, this formula allows us to determine the Wronskian associated with the solutions without effectively integrating the linear differential equation.

On the other hand, the Wronskian also occurs in the theory of space curves. For a smooth regular space curve (\ref{eq0}), the Wronskian $W(x,y,z)(t)$ of the functions defining parametrically the curve may arise, for instance, in the calculation of the distance $d(t)$ from the origin of the system of coordinates to the osculating plane at an arbitrary point of the curve while the Wronskian $W(x',y',z')(t)$ of their first derivatives (\ref{WronskianDer}) is used in the computation of the curve's torsion $\tau (t)$ as in (\ref{tau1}). Therefore, for example, any condition imposed upon on $d(t)$ and $\tau (t)$ yields a Wronskian-involving underdetermined nonlinear differential equation for the curve's defining functions.

In this paper, the class of nonlinear ODEs (\ref{EqnWronskian}) involving the Wronskians (\ref{Wronskian1}) and (\ref{WronskianDer}) of the unknown functions $x(t)$, $y(t)$, and $z(t)$ and, respectively, of their first derivatives $x'(t)$, $y'(t)$, and $z'(t)$ is introduced and analyzed in the context of the family of solutions satisfying the auxiliary linear differential equation (\ref{EqnU3}) such that the condition (\ref{ConditionWronsian}) holds. The motivation of considering differential equations of this kind lies in the study of the condition (\ref{condition}) that was
introduced for specific curves by the Romanian mathematician Gheorghe Tzitzeica in 1911 during his work on affine invariants \cite{Tzitzeica2}.
Nowadays, a space curve satisfying the relation (\ref{condition}) that may be written in an equivalent form as (\ref{TzitzeicaWronskian}) is called a \textit{Tzitzeica curve}. Along with Tzitzeica curves, Tzitzeica surfaces are affine invariants.
A \textit{Tzitzeica surface} is a surface for which the ratio of its Gaussian curvature and the fourth power of the distance from the origin to the tangent plane at any arbitrary point of the surface is constant \cite{Tzitzeica1}. It may be shown that the asymptotic curves on a Tzitzeica surface with negative Gaussian curvature are Tzitzeica curves
(see, e.g., \cite{Agnew}).
Although during the past years the Tzitzeica curves have been of interest to many geometers, their related nonlinear ODE has not been studied in detail so far, and, hence, there are known only a few examples of Tzitzeica curves defined explicitly in algebraic, transcendental, or integral forms (see \cite{Agnew}, \cite{BilaEniDSDG}, \cite{Crasmareanu}, and \cite{LWilliams_2013}). In \cite{Agnew}, Agnew et al. used the \texttt{Mathematica} software to find and illustrate the asymptotic curves on the Tzitzeica surface of revolution $z(x^2+y^2)=1$ and expressed them in cylindrical coordinates in terms of logarithmic and exponential functions. In \cite{Crasmareanu}, Cr\^{a}\c{s}m\u{a}reanu determined elliptic and hyperbolic cylindrical Tzitzeica curves written in integral form. In the paper by Williams \cite{LWilliams_2013}, the nonlinear Tzitzeica curve equation has been derived explicitly, and new closed-form solutions have been presented. B\^{\i}l\u{a} and Eni \cite{BilaEniDSDG} showed that the nonlinear ODE (\ref{TzitzeicaWronskian})
admits particular solutions obtained by augmenting it with a side condition consisting in a third-order homogeneous linear ODE with constant coefficients. In this paper, it is shown that the Tzitzeica curve equation also admits solutions that satisfy the side condition (\ref{EqnU3}). Although this is a slight generalization of the auxiliary condition introduced in \cite{BilaEniDSDG}, by
Proposition \ref{Prop1}, the new attached equation yields a larger family of solutions to the Tzitzeica curve equation. The new interesting solutions (\ref{sol3.1}) and (\ref{neweqn}) are expressed in closed-form or in terms of Airy's Ai and Bi functions. In Fig.~\ref{Figure1}, the software \texttt{Geogebra} has been used to visualize the Tzitzeica curve (\ref{neweqn}) on the Tzitzeica surface of equation $yz=1-4x^2$, surface that was found in \cite{Udriste_Bila_1999} by applying the classical Lie method to the Tzitzeica surface partial differential equation.
Additionally, for various functions $\gamma (t)$ in (\ref{EqnU3}), new Tzitzeica curves expressed in terms of other special functions \cite{Abramowitz} such as Bessel or generalized hypergeometric functions may be investigated. Theorem \ref{thm1} gives a generalization of these results for the class of nonlinear ODEs (\ref{EqnWronskian}).

The structure of the paper is the following. In Section \ref{Section2},
the side condition (\ref{EqnU3}) is explored in the case of the class of underdetermined nonlinear Wronskian-involving ODEs (\ref{EqnWronskian}) for which (\ref{ConditionWronsian}) holds. A new method for solving these types of nonlinear differential equations with the help of the auxiliary condition (\ref{EqnU3}) is introduced. The novel approach is exemplified in Section \ref{Section3} for the Tzitzeica curve equation (\ref{TzitzeicaWronskian}) for which intriguing solutions are obtained. In the last section, conclusions of this work are presented.

\section{A class of Wronskian-Involving Equations}\label{Section2}

\subsection{On a related side condition}\label{2.1}

If $x(t)$, $y(t)$, and $z(t)$ are smooth functions defined on an open interval $I\subset \mathbb{R}$, then their Wronskian is defined as
\begin{equation}\label{Wronskian1}
W(x,y,z)(t) = \left \vert
\begin{array}{ccc}
x(t) & y(t) & z(t) \\
x'(t)   & y'(t) & z'(t) \\
x''(t)  & y''(t) & z''(t)
\end{array}
\right \vert .
\end{equation}
Similarly, the Wronskian of the functions' first derivatives is given by
\begin{equation}\label{WronskianDer}
W(x',y',z')(t)= \left \vert
\begin{array}{ccc}
x'(t) & y'(t) & z'(t) \\
x''(t)   & y''(t) & z''(t) \\
x'''(t)  & y'''(t) & z'''(t)
\end{array}
\right \vert .
\end{equation}
In what follows, assume that
\begin{equation}\label{ConditionWronsian}
W(x,y,z)(t)\neq 0 \quad \text{and}\quad W(x',y',z')(t)\neq 0,
\end{equation}
for all $t\in I$. This supposition implies that the functions $x(t)$, $y(t)$, and $z(t)$ and, respectively, their first derivatives $x'(t)$, $y'(t)$, and $z'(t)$ are linearly independent (see \cite{BoyceDiPrima}, p. 221). Consider the differential equation
\begin{equation}\label{EqnWronskian}
{\cal {F}}\left (W(x,y,z)(t),W(x',y',z')(t)\right )=0
\end{equation}
involving the Wronskians (\ref{Wronskian1}) and (\ref{WronskianDer}), where ${\cal {F}}$ is a smooth real-valued function in two variables. The equation above is a third-order underdetermined differential equation in the unknown functions $x(t)$, $y(t)$, and $z(t)$. Next, the problem of finding solutions to (\ref{EqnWronskian}) is explored in the context of a specific side condition that would allow the Wronskian (\ref{WronskianDer}) to be expressed in terms of the Wronskian (\ref{Wronskian1}). Suppose that the  functions $x(t)$, $y(t)$, and $z(t)$ satisfy the following homogeneous third-order linear ODE
\begin{equation}\label{EqnU1}
u'''+\beta (t) u''+\gamma(t)u'+\delta  u=0,
\end{equation}
where $\beta (t)$ and $\gamma (t)$ are smooth functions defined on the interval $I$, and $\delta $ is a nonzero constant.
For simplicity, in this paper, $\delta $ is considered a nonzero constant. However, the method explained below may be modified to accommodate $\delta $ variable too.
Since the functions $x(t)$, $y(t)$, and $z(t)$ are linearly independent, they form a fundamental set of solutions to (\ref{EqnU1}). Therefore, the general solution to this equation may be written as $u(t)=C_1x(t)+C_2y(t)+C_3z(t)$, where $C_i$ are arbitrary real constants. If the leading derivatives
$x'''=-\beta (t) x''-\gamma (t)x'-\delta x$, $y'''=-\beta (t)y''-\gamma (t)y'-\delta y$, and $
z'''=-\beta (t)z''-\gamma (t)z'-\delta z$ are substituted
into (\ref{WronskianDer}), then after applying a few properties of determinants,
this Wronskian may be expressed as
\begin{equation}\label{WronskianDerNew}
W(x',y',z')(t)=-\delta W(x,y,z)(t), \quad t\in I.
\end{equation}
Replacing the above relation into the differential equation (\ref{EqnWronskian}) yields
\begin{equation}\label{EqnWronskian2}
{\cal {G}}\left (W(x,y,z)(t); \delta \right )=0,
\end{equation}
where ${\cal {G}}$ denotes the resulting function that depends on $W(x,y,z)(t)$ (the notation used here shows that  $\delta $ occurs as well in the reduced equation due to (\ref{WronskianDerNew})).
On the other hand, according to Abel-Liouville-Ostrogradski formula (see, e.g., \cite{BoyceDiPrima}, p. 239),
the Wronskian of the set of fundamental solutions of \eqref{EqnU1} satisfies the first-order linear ODE
\begin{equation}\label{N3}
W'(x,y,z)(t)=-\beta (t) W(x,y,z)(t).
\end{equation}
The differentiation of (\ref{EqnWronskian2}) with respect to $t$ leads to $W'(x,y,z)(t)=0$ for any $t\in I$. By (\ref{ConditionWronsian}), the Wronskian $W(x,y,z)(t)$ does not vanish, and, thus, the equation (\ref{N3}) implies $\beta (t) = 0$ for all $t\in I$.
In conclusion, the equation (\ref{EqnU1}) reduces to
\begin{equation}\label{EqnU3}
u'''+\gamma(t)u'+\delta  u=0,
\end{equation}
where $\gamma (t)$ is a smooth function and $\delta $ is a nonzero constant.
In this way, the equation (\ref{EqnWronskian}) has been reduced to the homogeneous third-order linear ODE (\ref{EqnU3}) along with the equation (\ref{EqnWronskian2}) provided that (\ref{ConditionWronsian}) holds. Observe that now the Wronskian $W(x,y,z)(t)=C_0$ is constant, and the equation (\ref{EqnWronskian2}) turns into the compatibility condition
\begin{equation}\label{EqnWronskian3}
{\cal {G}}\left (C_0; \delta \right )=0.
\end{equation}
In conclusion, the following result has been proven.

\begin{theorem}\label{thm1}
Any three linearly independent solutions $x(t)$, $y(t)$, and $z(t)$ of the homogeneous third-order linear ODE (\ref{EqnU3}) conditioned to (\ref{ConditionWronsian}) satisfy the nonlinear differential equation (\ref{EqnWronskian}) provided that (\ref{EqnWronskian3}) holds.
\end{theorem}

\subsection{Special side conditions}\label{2.2}

In what follows, the ODE (\ref{EqnU3}) is integrated in a few particular cases that yield to explicit or integral form solutions.
Recall that in each of these cases the Wronskian of the solutions is constant.
As it has been shown in the previous section, the equation (\ref{EqnWronskian}) has been reduced to the side condition (\ref{EqnU3}) and compatibility condition (\ref{EqnWronskian3}).

\begin{remark}\label{Exp1}
If $\gamma (t) =0$, then the equation (\ref{EqnU3}) becomes $u'''+\delta u=0$, and its set of fundamental solutions is
\begin{align}
x(t)&=\exp\left(-\delta ^{1/3} t\right), \nonumber \\
y(t)&=\exp\left(\frac{1}{2}\delta ^{1/3} t\right)\cos\left(\frac{\sqrt{3}}{2}\delta ^{1/3} t \right),\nonumber \\
z(t)&=\exp\left(\frac{1}{2}\delta ^{1/3} t\right)\sin \left(\frac{\sqrt{3}}{2}\delta ^{1/3} t \right), \label{sol1}
\end{align}
with $t \in \mathbb{R}$. In this case, $W(x,y,z)(t)= -3\delta \sqrt{3}/2$.
\end{remark}

\begin{remark}\label{Exp2}
For $\gamma (t)=\tilde {\gamma} $, where $\tilde {\gamma} $ is a nonzero real constant, the equation (\ref{EqnU3}) is reduced to the homogeneous linear ODE with constant coefficients
\begin{equation}
u'''+\tilde {\gamma} u'+\delta u=0 \label{case2}
\end{equation}
whose solutions have been discussed in detail in the context of the Tzitzeica curve equation in \cite{BilaEniDSDG}.
The characteristic equation related to (\ref{case2}) is the depressed cubic equation $v^3+\tilde {\gamma} v+\delta =0$.
The character of the solutions depends on the sign of the associated determinant
$$
D=-4\gamma ^3 - 27\delta ^2.
$$
\par
\noindent
a) If $D>0$, the depressed cubic equation has three real distinct solutions $v_i\neq 0$, $i=1,2,3$ whose sum is zero. The conditions $v_3\neq v_i$, with $i=1,2$, imply $v_2\neq -2v_1$ and $v_2\neq -v_1/2$.
In this case, the fundamental set of solutions of (\ref{case2}) is
\begin{equation}
x(t)=\exp\left(v_1 t\right), \quad
y(t)=\exp\left(v_2 t\right),\quad
z(t)=\exp\left(-\left(v_1+v_2\right) t\right),\label{sol2.2.1}
\end{equation}
with $t\in \mathbb{R}$, for which $W(x,y,z)(t)=(v_2-v_1)(2v_1+v_2)(v_1+2v_2)$.
\medskip
\par
\noindent
b) For $D=0$, the characteristic equation has three real solutions out of which two are equal, i.e., $v_1=v_2\neq 0$ and $v_3=-2v_1$. Therefore,
\begin{equation}
x(t)=\exp(v_1 t), \quad y(t)= t\exp(v_1 t), \quad z(t)=\exp(-2v_1 t), \label{sol2.2.2}
\end{equation}
with $t\in \mathbb{R}$, form a fundamental solution set to (\ref{case2}). In this case, their Wronskian is given by $W(x,y,z)(t)=9v_1^2$.
\medskip
\par
\noindent
c) If $D<0$, the depressed cubic equation has two nonreal complex conjugate solutions, $v_{1,2}=m \pm i n $ and one real nonzero solution
$v_3=-2m$, where $m, n\neq 0 $ are real numbers. The integration of (\ref{case2}) yields
\begin{equation}
x(t)=\exp(m t)\cos(n t),\;\;
y(t)= \exp(m t)\sin(n t), \;\;
z(t)= \exp(-2m t), \label{sol2.2.3}
\end{equation}
where $t\in \mathbb{R}$, and for which $W(x,y,z)(t)=n(9m^2+n^2)$.
\end{remark}

\begin{remark}\label{Exp3}
If $\gamma (t)=\delta t$, then (\ref{EqnU3}) may be integrated once and becomes $u''+\delta t u=C$, where $C\in \mathbb{R}$ (here $C$ is nonzero whenever three linearly independent solutions to (\ref{EqnU3}) are needed). The solutions of the latter ODE are given in integral form in terms of the Airy's Ai and Bi functions of the first and second kind, respectively (\cite{PolyaninZaitsev}, p. 214), that is,
\begin{align}
x(t)&=\text{Ai}\left (-\delta ^{1/3}t\right), \nonumber \\
y(t)&=\text{Bi}\left (-\delta ^{1/3}t\right), \nonumber \\
z(t)&= \pi \delta ^{-1/3}\left (x(t)\int y(s)ds -y(t)\int x(s)ds \right), \label{sol3.1}
\end{align}
with $t \in \mathbb{R}$, and their Wronskian is $W(x,y,z)(t)=-\delta ^{1/3}/ \pi $.
\end{remark}

\subsection{On an integration method involving a side condition}\label{2.3}

Assume that $x_0(t)$ is a known solution to ODE (\ref{EqnU3}). In this situation, a method for integrating the ODE (\ref{EqnU3}) may be introduced as follows.
\medskip
\par
\noindent
Step 1. Substitute $u=x_0(t)$ into the ODE (\ref{EqnU3}) and solve the equation for $\gamma (t)$. Denote $\gamma_0(t;\delta)$ the resulting solution that depends on $\delta $ too.
\medskip
\par
\noindent
Step 2. Replace $\gamma_0(t;\delta)$ in (\ref{EqnU3}) and obtain $u'''+\gamma_0(t;\delta )u'+\delta  u=0$. By using the change of function
\begin{equation}
u(t)=x_0(t)\int w(s)ds, \label{ChangeW}
\end{equation}
the previous equation is reduced to
\begin{equation}
x_0(t)w''+3x'_0(t)w'+\left[3x''_0(t)+\gamma_0(\delta;t)x_0(t)\right]w =0 \label{EqnW}
\end{equation}
which is a homogeneous linear second-order ODE for the function $w(t)$.
\medskip
\par
\noindent
Step 3. Integrate (\ref{EqnW}) and determine its general solution $w(t)=C_1w_1(t)+C_2w_2(t)$, where
$w_1(t)$ and $w_2(t)$ represent its fundamental solution set.
\medskip
\par
\noindent
Step 4. Replace $w(t)$ into (\ref{ChangeW}) to have
$$
u(t)=x_0(t)\int w(s) ds = C_1 x_0(t)\int w_1(s)ds + C_2 x_0(t)\int w_2(s)ds + C_3 x_0(t),
$$
where $t\in I$. In conclusion,
\begin{equation}
x(t)=x_0(t), \quad y(t)=x_0(t)\int w_1(s)ds, \quad z(t)=x_0(t)\int w_2(s)ds, \label{Soln1}
\end{equation}
with $t\in I$, form the fundamental set of solutions of the equation of (\ref{EqnU3}).

\section{Solutions to Tzitzeica Curve Equation}\label{Section3}

\subsection{Tzitzeica Curve Equation}\label{3.1}

In this section, the compatibility of the side condition (\ref{EqnU3}) for the underdetermined nonlinear ODE (\ref{EqnWronskian})
is presented in detail in the case of the Tzitzeica curve equation (\ref{TzitzeicaWronskian}). Consider
\begin{equation}
\mathbf{r}(t)=\left ( x(t), y(t), z(t)\right ),\quad t\in I  \label{eq0}
\end{equation}
a smooth, regular space curve, where $I\subset \mathbf{R}$ is an interval. Assume that the curve's curvature
$$
k(t) = \frac{\vert \vert  \mathbf{r'}(t)\times \mathbf{r''}(t)\vert \vert  }{\vert \vert  \mathbf{r'}(t)\vert \vert  ^3}, \quad t\in I
$$
and torsion
\begin{equation}
\tau (t)= \frac{\langle \mathbf{r'}(t),\mathbf{r''}(t),\mathbf{r'''}(t)\rangle }{\vert \vert  \mathbf{r'}(t)\times \mathbf{r''}(t)\vert \vert  ^2}, \quad t\in I    \label{tau1}
\end{equation}
do not vanish on $I$. Here $\vert \vert  \mathbf{r'}(t)\times \mathbf{r''}(t)\vert \vert $ is the magnitude of
the cross product of the tangent vector $\mathbf{r'}(t)$ and the acceleration vector $\mathbf{r''}(t)$,
 $\mathbf{r'}(t)\times \mathbf{r''}(t)$, $\vert \vert  \mathbf{r'}(t)\vert \vert $ is the magnitude of $\mathbf{r'}(t)$,
and
$$
\langle \mathbf{r'}(t),\mathbf{r''}(t),\mathbf{r'''}(t)\rangle  = \left \vert
\begin{array}{ccc}
x'(t) & y'(t) & z'(t) \\
x''(t)   & y''(t) & z''(t) \\
x'''(t)  & y'''(t) & z'''(t)
\end{array}
\right \vert
$$
is the mixed product of vectors $\mathbf{r'}(t)$, $\mathbf{r''}(t)$, and $\mathbf{r'''}(t)$ (see, for instance,
\cite{Presley}, p. 48). The curvature shows
the amount by which a curve deviates from being a line, and its torsion shows how sharply the curve is twisting out of its osculating plane.
If the torsion of the curve is nonzero then $W(x',y',z')(t)\neq 0$  and, hence, the functions $x'(t)$, $y'(t)$, and $z'(t)$ are linearly independent.
On the other hand, the osculating plane at an arbitrary point of the curve (\ref{eq0}) is the plane generated by the vectors $\mathbf{r'}(t)$ and $\mathbf{r''}(t)$. Therefore, its equation in the determinant form is given by
\[
\left \vert
\begin{array}{ccc}
x-x(t) & y-y(t) & z-z(t) \\
x'(t)   & y'(t) & z'(t) \\
x''(t)  & y''(t) & z''(t)
\end{array}
\right \vert  =0.
\]
If $d(t)$ denotes the distance from the origin to the osculating plane of the curve, it can be shown that
$$
d (t) = \frac{1}{\vert \vert  \mathbf{r'}(t) \times \mathbf{r''}(t)\vert \vert  }\left \vert
\begin{array}{ccc}
x(t) & y(t) & z(t) \\
x'(t)   & y'(t) & z'(t) \\
x''(t)  & y''(t) & z''(t)
\end{array}
\right \vert .
$$
The distance $d(t)$ does not vanish iff $W(x,y,z)(t)\neq 0$. In what follows, in addition to the condition that the curvature and the torsion of the curve (\ref{eq0}) are nonzero, the distance $d(t)$ is assumed nonzero on $I$ as well.  Therefore, the curve's defining functions $x(t)$, $y(t)$, and $z(t)$ are assumed to satisfy the condition (\ref{ConditionWronsian}).

A \textit{Tzitzeica curve} is defined as a space curve with the property that  the ratio of its torsion $\tau (t)$ and the square of the distance $d(t)$ from the origin to the osculating plane at an arbitrary point of the curve is constant, i.e,
\begin{equation}
\frac{\tau (t)}{d^2(t)}=\alpha , \label{condition}
\end{equation}
for all $t\in I$, where $\alpha \neq 0$ is a nonzero real number called here the \textit{curve's constant}. Substituting $\tau (t)$ and $d(t)$ into (\ref{condition}) yields
$$
\left \vert
\begin{array}{ccc}
x'(t) & y'(t) & z'(t) \\
x''(t)   & y''(t) & z''(t) \\
x'''(t)  & y'''(t) & z'''(t)
\end{array}
\right \vert
=\alpha \left \vert
\begin{array}{ccc}
x(t) & y(t) & z(t) \\
x'(t)   & y'(t) & z'(t) \\
x''(t)  & y''(t) & z''(t)
\end{array}
\right \vert ^2.
$$
The above equation can be rewritten in the terms of the Wronskians (\ref{Wronskian1}) and (\ref{WronskianDer}) as follows
\begin{equation}
W(x',y',z')(t)=\alpha \left [W(x,y,z)(t)\right ]^{2}. \label{TzitzeicaWronskian}
\end{equation}
Based on Theorem 2.1, the Wronskian-involving equation (\ref{TzitzeicaWronskian}) admits particular solutions satisfying the side condition (\ref{EqnU3}) provided that (\ref{ConditionWronsian}) holds.

Assume that the functions $x(t)$, $y(t)$, and $z(t)$ are solutions to the ODE (\ref{EqnU1}). In this case, the
Wronskian of their first derivatives (\ref{WronskianDer}) satisfies the relation (\ref{WronskianDerNew}) which may be used to rewrite the equation (\ref{TzitzeicaWronskian}) in the form
\begin{equation}
-\delta W(x,y,z)(t)=\alpha \left [W(x,y,z)(t)\right ]^{2}, \label{TzitzeicaCurveEqn}
\end{equation}
or, equivalently, as follows
$$
W(x,y,z)(t)\left[\alpha W(x,y,z)(t) -\delta )\right]=0.
$$
Thus, in this case, the above relation represents the compatibility condition (\ref{EqnWronskian3}) .
Taking into account that $W(x,y,z)(t)$ does not vanish on $I$, it follows
\begin{equation}
\alpha = -\frac{\delta }{W(x,y,z)(t)}. \label{eqn_alpha}
\end{equation}
Since the Tzitzeica curve equation (\ref{TzitzeicaCurveEqn}) is a particular case of the Wronskian-involving equation (\ref{EqnWronskian}), according to Theorem \ref{thm1}, the following result has been proven.

\begin{proposition}\label{Prop1}
Any three linearly independent solutions $x$, $y$, and $z$ of the homogeneous third-order linear ODE (\ref{EqnU3}) for which the condition (\ref{ConditionWronsian}) holds satisfy the Tzitzeica curve equation (\ref{TzitzeicaCurveEqn}), where the curve's constant is given by
(\ref{eqn_alpha}).
\end{proposition}

\subsection{Examples of Tzitzeica Curves}\label{3.2}

A few examples of Tzitzeica curves that derive from the side condition (\ref{EqnU3}) are presented below. However, more examples may be given by exploring other options for the function $\gamma (t)$ in this auxiliary condition.

\medskip
\par
\noindent
\textit{Example 1}.
According to Remark \ref{Exp1}, the functions (\ref{sol1}) represent a fundamental set of solutions for the side condition  (\ref{EqnU3}) in the case $\gamma (t)=0$. Since the Tzitzeica curve's equation (\ref{TzitzeicaCurveEqn}) reduces to (\ref{eqn_alpha}), then (\ref{sol1}) represents a Tzitzeica curve with $\alpha = -2\sqrt{3}/9$ that lies on the surface $S: x(y^2+z^2)=1$. After applying the affine transformation $T: \tilde x = z, \tilde y = x, \tilde z = y$ to the surface $S$, this becomes $\tilde S: z(x^2+y^2)=1$. It may be shown that $\tilde S$ is a Tzitzeica surface with negative Gaussian curvature, and, therefore, its asymptotic curves are Tzitzeica curves. Agnew et al \cite{Agnew} have found the asymptotic curves (expressed in cylindrical coordinates) for $\tilde S$  by using the software \texttt{Mathematica}. It may be shown that the transformed Tzitzeica curve (\ref{sol1}) via the transformation $T$ is an asymptotic curve on the surface $\tilde S$.

\medskip
\par
\noindent
\textit{Example 2}.
Consider the side condition (\ref{case2}) in Remark \ref{Exp2}. Then there are three classes of transcendental Tzitzeica curves that are classified in terms of the sign of the determinant $D$ of the depressed equation.
The solutions (\ref{sol2.2.1}), (\ref{sol2.2.2}), and (\ref{sol2.2.3}) represent Tzitzeica curves for which their associated curve's constants are given, respectively, by (\ref{eqn_alpha}). These curves have been found and discussed in detailed in \cite{BilaEniDSDG}.

\medskip
\par
\noindent
\textit{Example 3}.
The solution (\ref{sol3.1}) represents a new Tzitzeica curve for which, by (\ref{eqn_alpha}), the curve's constant is $\alpha = \pi \delta ^{2/3}$. Interestingly, this curve is expressed in terms of Airy Ai and Bi functions.
\medskip
\par
\begin{figure}[h]
\centering
\includegraphics[width=6in]{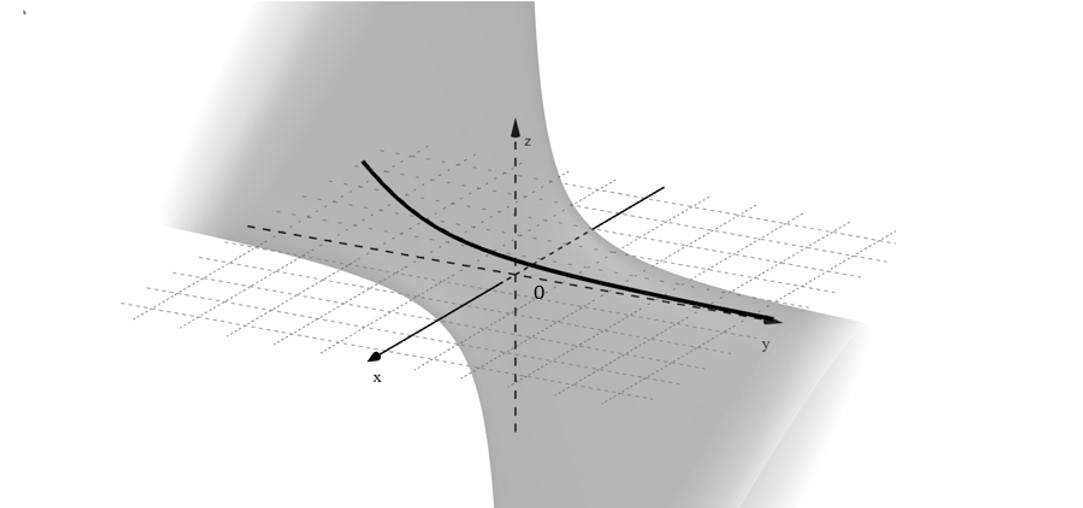}
\caption{\textit{The Tzitzeica curve defined by \eqref{neweqn} for $t\in [1.001,8]$ represented on the Tzitzeica surface of equation $yz=1-4x^2$.}}
\label{Figure1}
\end{figure}

\noindent
\textit{Example 4}. This example refers to the method introduced in Subsection \ref{2.3}. Assume that $x_0(t)=t^{-3/2}$ is a known solution to (\ref{EqnU3}). At Step 1, after substituting $x_0(t)$ back in the equation and solving for $\gamma (t)$, it follows
$$
\gamma (t)=\gamma _0(t;\delta )=\frac{8\delta t^3-105}{12t^2}.
$$
In particular, for $\delta = -27/8$, the ODE (\ref{EqnU3}) becomes
$$
u'''-\frac{9t^3+35}{4t^2}u'-\frac{27}{8}u=0.
$$
At Step 2, the change of variable (\ref{ChangeW}) reduces the above equation to
$$
w''-\frac{9}{2t}w'+\left(-\frac{9t}{4}+\frac{5}{2t^2}\right)w=0.
$$
The set of fundamental solutions is found at Step 3 to be
$$
w_1(t)=\frac{3}{2}\left(t^2-\sqrt{t}\right)\exp\left({t^{3/2}}\right),\quad
w_2(t)=-\frac{3}{2}\left(t^2+\sqrt{t}\right)\exp\left(-{t^{3/2}}\right),
$$
for $t>1$. Therefore, at Step 4, the following Tzitzeica curve is found
\begin{align}
x(t)&=t^{-3/2},\nonumber \\
y(t)&=\left(1-2t^{-3/2}\right)\exp\left({t^{3/2}}\right), \nonumber \\
z(t)&=\left(1+2t^{-3/2}\right)\exp\left({-t^{3/2}}\right), \label{neweqn}
\end{align}
where $t>1$. The Wronskian (\ref{Wronskian1}) becomes $W(x,y,z)(t)=27/4$, and, from (\ref{eqn_alpha}), the curve's constant takes the value $\alpha = 1/2$. The new Tzitzeica curve (\ref{neweqn}) lies on the Tzitzeica surface
$M:\;yz=1-4x^2$ with negative Gaussian curvature, surface that was found by using the symmetry analysis theory in \cite{Udriste_Bila_1999} (see Fig.~\ref{Figure1}). The curve (\ref{neweqn}) is not an asymptotic curve on $M$ as this is an hyperboloid with one sheet whose asymptotic curves are lines.

\section{Discussion and Conclusions}\label{section4}

The motivation of this work is based on the study of the Tzitzeica curve equation (\ref{TzitzeicaCurveEqn}) which is an underdetermined nonlinear autonomous ordinary differential equation arising from the condition stating that the ratio of the curve's torsion $\tau (t)$ and
the square of the distance $d (t)$ from the origin to its osculating plane at an arbitrary point of the curve is constant.
This paper is a continuation of author's work on closed-form solutions for the Tzitzeica curves that she has initiated with her students (see \cite{BilaEniDSDG} and \cite{LWilliams_2013}). It is shown that the auxiliary equation (\ref{EqnU3}) provides a large family of solutions for the class of Wronskian-involving equations (\ref{EqnWronskian}) which includes the Tzitzeica curve nonlinear ODE (\ref{TzitzeicaCurveEqn}). The side condition (\ref{EqnU3}) has been derived by using a linear combination of derivatives of the unknown functions that would allow specific properties of determinants to be applied. Another
key observation is the relation (\ref{WronskianDerNew}) that shows that the Wronskian (\ref{WronskianDer}) may be expressed in terms of the Wronskian (\ref{Wronskian1}). In Section \ref{Section2}, a method based on the side condition (\ref{EqnU3}) is introduced and allows one to find a Tzitzeica curve if a solution $x_0(t)$ to (\ref{EqnU3}) is known. As a consequence, new solutions for (\ref{EqnWronskian}) and, in particular, for (\ref{TzitzeicaCurveEqn}) are obtained.  It is intriguing how the auxiliary condition (\ref{EqnU3}) may provide a large variety of solutions to the nonlinear ODEs (\ref{EqnWronskian}) and, hence, to (\ref{TzitzeicaCurveEqn}). In the future, it will be interesting to make the connection of these results with the classical Lie symmetries associated with the Tzitzeica curve equation and with other particular Wronskian-type equations.

\textbf{Acknowledgments}. Part of these results have been presented to the online XVth International Conference on Differential Geometry and Dynamical Systems, held on August 26--29, 2021, Bucharest, Romania, and, hence, the author would like to thank the organizers, in particular, Prof. Dr. C. Udri\c{s}te and Prof. Dr. V. B\u{a}lan. The author would also like to thank Prof. Dr. M. Cr\^{a}\c{s}m\u{a}reanu for reading the manuscript and making interesting comments that have improved the presentation of the paper.


\begin{thebibliography}{20}


\bibitem{Abramowitz}
Abramowitz, M. and Stegun, I., eds. Handbook of Mathematical Functions. New York: Dover (1972)

\bibitem{Agnew}
A. F. Agnew, A. Bobe, W. G. Boskoff and B. D. Suceava, {\em Tzitzeica curves and surfaces}, The Mathematica Journal, \textbf{12}, 1-18 (2010)


\bibitem{BilaEniDSDG}
N. B\^{\i}l\u{a} and M. Eni, {\em Particular solutions to the Tzitzeica curve equation}, Differential Geometry - Dynamical Systems, \textbf{24}, 38-47 (2022)

\bibitem{BoyceDiPrima}
W. E. Boyce and R. C. DiPrima, {\em Elementary Differential Equations and Boundary Value Problems}, 8th edition, John Wiley \& Sons, Inc. (1986)

\bibitem{Crasmareanu}
M. Cr\^{a}\c{s}m\u{a}reanu, {\em Cylindrical Tzitzeica curves implies forced harmonic oscillators}, Balkan J.
Geom. Appl.,\textbf{7}, No. 1, 37-42 (2002)

\bibitem{Muir}
T. Muir, {\em A treatise on the theorie of determinants}, London. Macmillan (1882)

\bibitem{OlverBook}
P. J. Olver, {\em Applications of Lie Groups to Differential Equations},
Graduate Texts in Mathematics, vol. 107, Springer-Verlag, New York (1986)

\bibitem{PolyaninZaitsev}
A. D. Polyanin, V. F. Zaitsev, {\em
Handbook of Exact Solutions for Ordinary Differential Equations} Second Edition, Chapman \& Hall/CRC Press, Boca Raton (2003)

\bibitem{Presley}
A. Pressley, {\em Elementary Differential Geometry}, Springer Undergraduate Mathematics Series, Springer-Verlag London Limited
(2012)

\bibitem{Tzitzeica1}
G. Tzitzeica, {\em Sur une nouvelle classes de surfaces}, Comptes Rendus, Acad. Sci. Paris, Paris, \textbf{144}, 1257-1259 (1907)


\bibitem{Tzitzeica2}
G. Tzitzeica, {\em Sur certaines courbes gauches}, Ann. de l'Ec. Normale Sup., \textbf{28}, 9-32 (1911)

\bibitem{Udriste_Bila_1999}
C. Udri\c{s}te, N. B\^{\i}l\u{a}, {\em Symmetry group of Tzitzeica surfaces PDE}, Balkan J. Geom.
Appl., \textbf{4}(2), 123-140 (1999)


\bibitem{LWilliams_2013}
L. R. Williams, {\em On The Tzitzeica Curve Equation}, Explorations: The Undergraduate Research and Creative Activities Journal for the State of North Carolina, \textbf{VIII}, 105-115 (2013)

\end{thebibliography}
\end{document}